\newenvironment{E}{\begin{equation}}{\end{equation}}
\def\proof{\noindent{\bf Proof: }}
\def\qed{ \hskip 20pt{\vrule height7pt width6pt depth0pt}\hfil}
\def\forb{{\hbox{forb}}}
\def\ex{{\hbox{ex}}}
\def\supp{{\hbox{supp}}}
\def\even{{\hbox{even}}}
\def\odd{{\hbox{odd}}}
\def\inc{{\hbox{Inc}}}

\def\cF{{\mathcal{F}}}

\def\0{{\bf 0}}
\def\1{{\bf 1}}

\def\Av{{\mathrm{Avoid}}}
\newcommand{\linelessfrac}[2]{\genfrac{}{}{0pt}{}{#1}{#2}}
\newcommand{\ncols}[1]{\| #1 \|}
\newcommand{\rf}[1]{(\ref{#1})}
\newcommand{\trf}[1]{Theorem~\ref{#1}}
\newcommand{\lrf}[1]{Lemma~\ref{#1}}

\newcommand{\corf}[1]{Conjecture~\ref{#1}}

\newcommand{\srf}[1]{Section~\ref{#1}}
\newcommand{\remrf}[1]{Remark~\ref{#1}}
\newcommand{\tbrf}[1]{Table~\ref{#1}}
\documentclass[12pt]{article}
\newtheorem{thm}{Theorem}[section]
\newtheorem{lemma}[thm]{Lemma}

\newtheorem{cor}[thm]{Corollary}
\newtheorem{conj}[thm]{Conjecture}

\newtheorem{remark}[thm]{Remark}

\renewcommand{\l}{\ell}
\usepackage{amssymb}
\usepackage{amsmath}
\title{ Forbidden Families of Configurations}
\author{R.P. Anstee\thanks{Research supported in part by NSERC}, Christina L. Koch\thanks{Research suppoorted in part by NSERC of first author} \\Mathematics Department\\The University of British Columbia\\Vancouver, B.C. Canada V6T 1Z2\\ {\small\texttt{anstee@math.ubc.ca}} {\small\texttt{clkoch@math.ubc.ca}}
\\\mbox{\ }}

\begin{document}
\maketitle
\begin{abstract}
A \emph{simple} matrix is a (0,1)-matrix with no repeated columns. For a (0,1)-matrix $F$, we say
 that a (0,1)-matrix $A$ has $F$ as a \emph{configuration} if there is a submatrix of $A$ which is a row and column permutation of $F$ (\emph{trace} is the set system version of a \emph{configuration}). Let $\ncols{A}$ denote the number of columns of $A$. Let ${\cal F}$ be a family of matrices. We define the extremal function $\forb(m,{\cal F})=\max\{\ncols{A}\,:\, A \hbox{ is }m\hbox{-rowed simple matrix and has no configuration }F\in{\cal F}\}$.
We consider some families ${\cal F}=\{F_1,F_2,\ldots ,F_t\}$   such that individually each $\forb(m,F_i)$ has greater asymptotic growth than  $\forb(m,{\cal F})$.
\end{abstract}

\section{Introduction}

We are initiating an exploration of families of forbidden configurations in this paper as recommended in \cite{miguelthesis}. We need some notation. 
We say a matrix is \emph{simple} if it is a (0,1)-matrix with no repeated columns. Such a matrix can be viewed as an element-set incidence matrix.
Given two (0,1)-matrices $F,A$, if there is a submatrix of $A$ which is a row and column permutation of $F$ then we say $A$ has $F$ as a \emph{configuration} and write $F\prec A$.
In set terminology we could use the notation \emph{trace}. For a subset of rows $S$, we define $A|_S$ as the submatrix of $A$ consisting of rows $S$ of $A$.  We define $[n]=\{1,2,\ldots ,n\}$. If $F$ has $k$ rows and $A$ has $m$ rows and $F\prec A$ then there is a $k$-subset $S\subseteq[m]$ such that $F\prec A|_S$. For two $m$-rowed matrices $A,B$ we use $[A\,|\,B]$ to denote the concatenation of $A,B$ yielding a larger $m$-rowed matrix. We define $t\cdot A$ as the matrix obtained from concatenating $t$ copies of $A$. These two operations need not yield simple matrices. Let $A^c$ denote the (0,1)-complement of $A$.

Define $\ncols{A}$ as the number of columns of $A$. For some set of matrices  ${\cal F}$, we define our extremal problem as follows:
$$\Av(m,{\cal F})=\{A\,:\,A\hbox{ is $m$-rowed, simple,}F\not\prec A\hbox{ for all }F\in{\cal F}\},$$ 
$$\forb(m,{\cal F})=\max_A\{\ncols{A}\,:\,A\in\Av(m,{\cal F})\}.$$
When $|{\cal F}|=1$ and ${\cal F}=\{F\}$, we write $\Av(m,F)$ and $\forb(m,F)$.
  A conjecture of Anstee and Sali \cite{AS05} for a single configuration sometimes makes the correct predictions for the asymptotic growth of $\forb(m,{\cal F})$.  Let $I_k$ denote the $k\times k$ identity matrix and let $T_k$ denote the $k\times k$ triangular simple matrix  with a 1 in position $(i,j)$ if and only $i\le j$. 
 For an $m_1\times n_1$ simple matrix $A$ and a $m_2\times n_2$ simple matrix $B$, we define the 2-fold product $A\times B$ as the $(m_1+m_2)\times n_1n_2$ simple matrix whose columns are obtained from placing a column of $A$ on top of a column of $B$ in all possible ways. This generalizes to $p$-fold products.
For a configuration $F$ we define $X(F)$ as the smallest value of $p$ such that 
$F\prec A_1\times A_2\times\cdots\times A_p$ for  every $p$-fold product where $A_i\in\{I_{m/p},I_{m/p}^c,T_{m/p}\}$.

\begin{conj}\label{grand}\cite{AS05} We believe that $\forb(m,F)$ is $\Theta(m^{X(F)-1})$. \end{conj}

We think that the conjecture will help in guessing asymptotic  bounds for $\forb(m,{\cal F})$.
We may define $X({\cal F})$ as the smallest value of $p$ such that for every every $p$-fold product
$A_1\times A_2\times\cdots\times A_p$  where $A_i\in\{I_{m/p},I_{m/p}^c,T_{m/p}\}$ we have some $F\in{\cal F}$ with $F\prec A_1\times A_2\times\cdots\times A_p$.

Two easy remarks are the following. 

\begin{remark}We have that $\forb(m,\{F_1,F_2,\ldots,F_t\})=\forb(m,\{F_1^c,F_2^c,\ldots,F_t^c\})$.\label{complement}\end{remark}

\begin{remark}\label{FsubsetG}Let ${\cal F}\subseteq{\cal G}$. Then $\forb(m,{\cal G})\le \forb(m,{\cal F})$.\end{remark}

\begin{remark}\label{minimal}Let ${\cal F}$ be given with $F\in{\cal F}$. Let $F'$ be given with $F\prec F'$, Then  $ \forb(m,{\cal F}\cup\{F'\})=\forb(m,{\cal F})$.\end{remark}

In view of \remrf{minimal}, we define ${\cal F}$ to be \emph{minimal} if there are no pair $F,F'\in{\cal F}$ with $F\prec F'$.

Some examples are in order.
 Balanced and totally balanced matrices are classes of matrices which can each be defined using an infinite family of forbidden configurations.
Let $C_k$ denote the vertex-edge incidence matrix of the cyle of length $k$. Thus
$$\hbox{e.g. }\,\,C_3=\left[\begin{array}{@{}ccc@{}}
1&0&1\\ 
1&1&0\\
0&1&1\\   \end{array}\right], C_4=\left[\begin{array}{@{}cccc@{}}
1&0&0&1\\ 
1&1&0&0\\
0&1&1&0\\
0&0&1&1\\   \end{array}\right].$$
A  matrix $A$  is \emph{balanced} if has no configuration $C_k$ for $k$ odd and a matrix is \emph{totally balanced} if it has no configuration $C_k$ for all $k\ge 3$.
These are important classes of matrices. While the definitions do not require the matrices to be simple, it is still of interest how many different columns can there be in a balanced (resp. totally balanced) matrix on $m$ rows.  We obtain an upper bound using \remrf{FsubsetG} and the lower bound  follows from the result that any $m\times \forb(m,C_3)$ matrix $A\in\Av(m,C_3)$ is necessarily  totally balanced.

\begin{thm} \cite{A80}
 We have that $\forb(m,C_3)=\forb(m,\{C_3,C_4,C_5,C_6,\ldots\})=$\hfil\break $\forb(m,\{C_3,C_5,C_7,C_9,\ldots\})$.\end{thm} 
 
The result $\forb(m,C_3)=\binom{m}{2}+\binom{m}{1}+\binom{m}{0}$ is due to Ryser \cite{ryser72}. We note that \hfil\break $X(\{C_3,C_4,C_5,C_6,\ldots\})=X(\{C_3,C_5,C_7,C_9,\ldots\})=3$ where the construction $T_{m/2}\times T_{m/2}$ avoids $C_k$ for all $k\ge 3$. From another point of view, the result suggests that the bound for   a forbidden family  might arise from the most restrictive configuration in the family (i.e. $\forb(m,{\cal F})=\min_{F\in{\cal F}}\forb(m,F)$ or its asymptotic equivalent)  but this is generally not true.  The following examples suggest that forbidden families can behave  quite differently.  We consider the fundamental extremal function $\ex(m,H)$ which denotes the maximum number of edges in a (simple) graph on $m$ vertices that has no subgraph $H$.  Let $\1_k$ denote the $k\times 1$ column of 1's. We can connect this to forbidden families as follows. We note that $A\in\Av(m,\1_3)$ consists of at most $m+1$ the columns of column sum 0 or 1 and $A$ may have columns of sum 2. The columns of sum 2  can be interpreted as a vertex-edge incidence matrix  of a graph.  For a graph $H$, let $\inc(H)$ denote its vertex-edge incidence matrix. We deduce the following.
 
 \begin{lemma}We have that $\forb(m,\{\1_3,\inc(H)\})=\ex(m,H)+m+1$.\label{extendgraph}\end{lemma}
 
 Two sample results concerning $\ex(m,H)$ yield the following where the vertex-edge incidence matrix of the cycle of length $k$ is $C_k$.
 
 \begin{thm}\cite{KST} We have that $\forb(m,\{\1_3,C_4\})$ is $\Theta(m^{3/2})$.\label{C4}\end{thm}
  
  \begin{thm}\cite{BS} We have that $\forb(m,\{\1_3,C_6\})$ is $\Theta(m^{4/3})$.\label{C6}\end{thm} 
  
   Simonovits refers to an unpublished upper bound of Erd\H{o}s as the `Even Circuit Theorem' so the origins of the results are partly folklore.
   \corf{grand} is failing spectacularly on these examples ($X(\{\1_3,C_4\})=X(\{\1_3,C_6\})=2$) and also on the following example. You might note that $I_2\times I_2$ is the same as $C_4$ after a row and column permutation.
  
  \begin{thm}\cite{AKRS} We have that $\forb(m,\{I_2\times I_2, T_2\times T_2\})$ is $\Theta(m^{3/2})$.\label{unusual}\end{thm}

Balogh and Bollob\'as proved the following useful bound which is consistent with \corf{grand}. For fixed $k$, we have $X(\{I_k,I_k^c,T_k\})=1$ since all 1-fold products contain some element of $\{I_k,I_k^c,T_k\}$.
. 
\begin{thm} \cite{BB} Let $k$ be given. Then there is a constant {$c_k$} so that
$\forb(m,\{I_k,I_k^c,T_k\})={c_k}$.\label{constant}\end{thm}

 The following lemma is straightforward and  quite useful. 

\begin{lemma} Let ${\cal F}=\{F_1,F_2,\ldots ,F_k\}$ and ${\cal G}=\{G_1,G_2,\ldots , G_{\ell}\}$. Assume  that for every $G_i$, there is some $F_j$ with $F_j\prec G_i$.
Then $\forb(m,{\cal F})\le \forb(m,{\cal G})$.
\label{FinG}\end{lemma}

\proof Assume $\ncols{A}>\forb(m,{\cal G})$. Then for some $i\in [t]$, $G_i\prec A$. But by hypothesis there is some $F_j\in{\cal F}$ with $F_j\prec G_i$. But then 
$F_i\prec A$, verifying that $\forb(m,{\cal F})\le \forb(m,{\cal G})$.\qed
\vskip 5pt
Now combining with \trf{constant}, we obtain a surprising  classification.

\begin{thm} Let ${\cal F}=\{F_1,F_2,\ldots ,F_t\}$ be given. There are two possibilities. Either $\forb(m,{\cal F})$ is $\Omega(m)$ or there exist 
$\ell,i,j,k$  with $F_i\prec I_{\ell}$,  with $F_j\prec I_{\ell}^c$ and  with $F_k\prec T_{\ell}$ in which case  there is a constant $c$ with $\forb(m,{\cal F})= c$.\label{classify}\end{thm}

\proof  Let $F_i$ be $a_i\times {b}_i$ and let $\ell=\max_{i\in [t]}(a_i+{b}_i)$. Let ${\cal G}=\{I_{\ell},I_{\ell}^c,T_{\ell}\}$. Then $F_j\nprec I_{\ell}$ implies 
$F_j\not\prec I_m$ for any $m\ge {\ell}$. Thus if $F_j\nprec I_{\ell}$ for $j=1,2,\ldots t$, then $\forb(m,{\cal F})$ is $\Omega(m)$ using the construction  $I_m$. The same holds for $I^c$ and $T$.
\qed

\vskip 10pt 
This paper considers all pairs of forbidden configurations drawn from \tbrf{tab:minquad}. The listed nine configurations are  \emph{minimal quadratic} configurations, namely those $Q$ for which $\forb(m,Q)$ is $\Theta(m^2)$ yet for any submatrix $Q'$ of $Q$, where $Q'\ne Q$, has $\forb(m,Q')$ being $O(m)$.   The minimal quadratic configurations of \tbrf{tab:minquad} have the virtue of having few possible 2-fold constructions avoiding them and so avoiding the configurations in pairs (or larger families) results in interesting interactions. \tbrf{tab:minquad} lists all the product constructions that yield the quadratic lower bounds which you can use this to compute $X(\{Q_i,Q_j\})$ for pairs $Q_i,Q_j$ in the table. The asymptotic growth rates of $\forb(m,\{Q_i,Q_j\})$ are collected together in \tbrf{tab:bounds} and the complete  analysis for any non-empty ${\cal F}\subset\{Q_1,Q_2,\ldots ,Q_9\}$ is in \trf{anysubset}. \srf{mainbds} handles those pairs with $X(\{Q_i,Q_j\})=3$ for which it is immediate that $\forb(m,\{Q_i,Q_j\})$ is $\Theta(m^2)$. Also we consider those cases where \lrf{FinG} when applied with \trf{constant} yield that $\forb(m,\{Q_i,Q_j\})$ is $O(1)$.
\srf{graphtheory} considers how to apply \lrf{extendgraph} more generally to help with $\forb(m,\{Q_5,Q_j\})$. \srf{SI} provides a new standard induction introduced in \cite{AL} that is useful in this context and helps with $\forb(m,\{Q_8,Q_j\})$ and  $\forb(m,\{Q_3,Q_j\})$. \srf{structure} considers the structures that arise from forbidding $Q_9$ and then uses this to obtain results on 
$\forb(m,\{Q_9,Q_j\})$.

\begin{table}
\begin{tabular}[c]{|c|c|c|p{3cm}|p{2cm}|}
  \hline
   & Configuration $Q_i$ &  $\forb(m,Q_i)$ & Construction(s) & Reference \\ \hline
  $Q_1$ & $\left[
             \begin{array}{cc}
               0 & 0 \\
               0 & 0 \\
             \end{array}
           \right]
  $ &   $\binom{m}{2} + \binom{m}{1} + \binom{m}{0}$ & $I^c \times I^c$ & \cite{G} \\ \hline
  $Q_2$ & $\left[\begin{array}{cc}
                                             1 & 1 \\
                                             1 & 1
                                           \end{array}
  \right]$ &  $\binom{m}{2} + \binom{m}{1} + \binom{m}{0}$ & $I \times I$ & \cite{G} \\ \hline
  $Q_3$ & $\left[\begin{array}{cccccc}
                       0 & 0 & 0 & 1 & 1 & 1 \\
                       0 & 1 & 1 & 0 & 0 & 1
                     \end{array}
  \right]$ &  $\lfloor\frac{m^2}{4}\rfloor + m + 1$ & $I \times I^c$ & \cite{AFS01} \\ \hline
  $Q_4$ & $\left[ \begin{array}{c}
             0 \\
             0 \\
             0
           \end{array} \right]
  $ &  $\binom{m}{2} + \binom{m}{1} + \binom{m}{0}$ & $I^c \times I^c$ & \cite{Sauer, PS, VC} \\ \hline
  $Q_5$ & $\left[ \begin{array}{c} 1 \\ 1 \\ 1 \end{array}\right]$ &  $\binom{m}{2} + \binom{m}{1} + \binom{m}{0}$ & $I \times I$ & \cite{Sauer, PS, VC} \\ \hline
  $Q_6$ & $\left[\begin{array}{ccc}
           1 & 0 & 0 \\
           0 & 1 & 0 \\
           0 & 0 & 1
         \end{array}
 \right]$ & $\binom{m}{2} + \binom{m}{1} + \binom{m}{0}$ & $\begin{matrix}I^c \times I^c\\ I^c \times T\\ T \times T\\ \end{matrix}$ & \cite{ryser72} \\ \hline
  $Q_7$ & $\left[\begin{array}{ccc}
                        0 & 1 & 1 \\
                        1 & 0 & 1 \\
                        1 & 1 & 0
                      \end{array}
  \right]$ &  $\binom{m}{2} + \binom{m}{1} + \binom{m}{0}$ & $\begin{matrix}I \times I\\ I \times T\\ T \times T\\ \end{matrix}$ & \cite{ryser72} \\ \hline
  $Q_8$ & $\left[\begin{array}{cccc}
                        1 & 0 & 1 & 0 \\
                        0 & 1 & 0 & 1 \\
                        0 & 0 & 1 & 1
                      \end{array}
   \right]$ &  $\lfloor \frac{m^2}{4} \rfloor + m + 1$ & $T \times T$  & \cite{survey} \\ \hline
  $Q_9$ & $\left[\begin{array}{cc} 1 & 0 \\ 1& 0 \\ 0 & 1 \\ 0 & 1 \end{array} \right]$ &
   $\binom{m}{2} + 2m - 1$ & $I \times T$ \newline $I^c \times T$  & \cite{FFP87}  \\
  \hline
\end{tabular}
\caption{Minimal Quadratic Configurations}
\label{tab:minquad}
\end{table}

 \begin{table}
\begin{tabular}[c]{|c|c|c|c|c|c|c|c|c|}
\hline
&$Q_2$&$Q_3$&$Q_4$&$Q_5$&$Q_6$&$Q_7$&$Q_8$&$Q_9$\\
\hline
$Q_1$&$\begin{matrix}\Theta(1)\\ \hbox{Th}~\ref{constantbd}\\ \end{matrix}$&$\begin{matrix}\Theta(m)\\ \hbox{Cor~\ref{orderbd}}\\ \end{matrix}$
&$\begin{matrix}\Theta(m^2)\\ \hbox{Th}~\ref{IxI}\\ \end{matrix}$
&$\begin{matrix}\Theta(1)\\ \hbox{Th}~\ref{constantbd}\\ \end{matrix}$&$\begin{matrix}\Theta(m^2)\\ \hbox{Th~\ref{IxI}}\\ \end{matrix}$&
$\begin{matrix}\Theta(1)\\ \hbox{Th}~\ref{constantbd}\\ \end{matrix}$
&$\begin{matrix}\Theta(m)\\ \hbox{Cor~\ref{F1F8,F2F8}}\\ \end{matrix}$&$\begin{matrix}\Theta(m)\\ \hbox{Cor~\ref{Q9}}\\ \end{matrix}$\\
\hline
$Q_2$&&$\begin{matrix}\Theta(m)\\ \hbox{Cor~\ref{orderbd}}\\ \end{matrix}$&$\begin{matrix}\Theta(1)\\ \hbox{Th}~\ref{constantbd}\\ \end{matrix}$
&$\begin{matrix}\Theta(m^2)\\ \hbox{Th}~\ref{IxI}\\ \end{matrix}$
&$\begin{matrix}\Theta(1)\\ \hbox{Th}~\ref{constantbd}\\ \end{matrix}$&$\begin{matrix}\Theta(m^2)\\ \hbox{Th}~\ref{IxI}\\ \end{matrix}$
&$\begin{matrix}\Theta(m)\\ \hbox{Cor~\ref{F1F8,F2F8}}\\ \end{matrix}$&$\begin{matrix}\Theta(m)\\ \hbox{Cor~\ref{Q9}}\\ \end{matrix}$\\
\hline
$Q_3$&&&$\begin{matrix}\Theta(m)\\ \hbox{Th~\ref{graphbd}}\\ \end{matrix}$&$\begin{matrix}\Theta(m)\\ \hbox{Th~\ref{graphbd}}\\ \end{matrix}$
&$\begin{matrix}\Theta(m)\\ \hbox{Cor~\ref{orderbd}}\\ \end{matrix}$&$\begin{matrix}\Theta(m)\\ \hbox{Cor~\ref{orderbd}}\\ \end{matrix}$
&$\begin{matrix}\Theta(m)\\ \hbox{Cor~\ref{Q3Q8}}\\ \end{matrix}$&$\begin{matrix}\Theta(m)\\ \hbox{Cor~\ref{Q9}}\\ \end{matrix}$\\
\hline
$Q_4$&&&&$\begin{matrix}\Theta(1)\\ \hbox{Th}~\ref{constantbd}\\ \end{matrix}$&$\begin{matrix}\Theta(m^2)\\ \hbox{Th}~\ref{IxI}\\ \end{matrix}$
&$\begin{matrix}\Theta(1)\\ \hbox{Th}~\ref{constantbd}\\ \end{matrix}$&$\begin{matrix}\Theta(m)\\ \hbox{Th~\ref{graphbd}}\\ \end{matrix}$
&$\begin{matrix}\Theta(m)\\ \hbox{Th~\ref{graphbd}}\\ \end{matrix}$\\
\hline
$Q_5$&&&&&$\begin{matrix}\Theta(1)\\ \hbox{Th}~\ref{constantbd}\\ \end{matrix}$&$\begin{matrix}\Theta(m^2)\\ \hbox{Th}~\ref{IxI}\\ \end{matrix}$
&$\begin{matrix}\Theta(m)\\ \hbox{Th~\ref{graphbd}}\\ \end{matrix}$&$\begin{matrix}\Theta(m)\\ \hbox{Th~\ref{graphbd}}\\ \end{matrix}$\\
\hline
$Q_6$&&&&&&$\begin{matrix}\Theta(m^2)\\ \hbox{Th}~\ref{TxT}\\ \end{matrix}$&$\begin{matrix}\Theta(m^2)\\ \hbox{Th}~\ref{TxT}\\ \end{matrix}$
&$\begin{matrix}\Theta(m^2)\\ \hbox{Th}~\ref{IxT}\\ \end{matrix}$\\
\hline
$Q_7$&&&&&&&$\begin{matrix}\Theta(m^2)\\ \hbox{Th}~\ref{TxT}\\ \end{matrix}$&$\begin{matrix}\Theta(m^2)\\ \hbox{Th}~\ref{IxT}\\ \end{matrix}$\\
\hline
$Q_8$&&&&&&&&$\begin{matrix}\Theta(m)\\ \hbox{Th~\ref{Q8Q9}}\\ \end{matrix}$\\
\hline

\end{tabular}
\caption{Asymptotic growth rates of $\forb(m,\{Q_i,Q_j\})$.}
\label{tab:bounds}
\end{table}

\section{Quadratic and Constant Bounds}\label{mainbds}

 First we are interested in pairs with $X(\{Q_i,Q_j\})=3$ for which it follows that 
 \hfil\break$\forb(m,\{Q_i,Q_j\})$ is $\Theta(m^2)$ (the upper bound follows from \lrf{FsubsetG} using that
  \hfil\break$\forb(m,\{Q_i\})$ is $O(m^2)$ for all $i\in [9]$).

\begin{thm} We have that $\forb(m,\{Q_1,Q_4,Q_6\})=\forb(m,\{Q_2,Q_5,Q_7\})$ is $\Theta(m^2)$.\label{IxI}\end{thm}
\proof We use the construction $I_{m/2}^c\times I_{m/2}^c\in\Av(m,\{Q_1,Q_4,Q_6\})$ to deduce that
$X(\{Q_1,Q_4,Q_6\})=3$ and $I_{m/2}\times I_{m/2}\in\Av(m,\{Q_2,Q_5,Q_7\})$ 
yields
$X(\{Q_2,Q_5,Q_7\})=3$.\qed

\begin{thm} We have that $\forb(m,\{Q_6,Q_7,Q_8\})$ is $\Theta(m^2)$.\label{TxT}\end{thm}
\proof The construction $T_{m/2}\times T_{m/2}\in\Av(m,\{Q_6,Q_7,Q_8\})$ shows that \hfil\break
$X(\{Q_6,Q_7,Q_8\})=3$.\qed

\begin{thm} We have that $\forb(m,\{Q_6,Q_9\})$ and $\forb(m,\{Q_7,Q_9\})$ are $\Theta(m^2)$.\label{IxT}\end{thm}
\proof We use the construction $I_{m/2}^c\times T_{m/2}\in \Av(m,\{Q_6,Q_9\})$ 
to deduce that
$X(\{Q_6,Q_9\})=3$
and $I_{m/2}\times T_{m/2}\in \Av(m,\{Q_7,Q_9\})$ yields $X(\{Q_7,Q_9\})=3$.\qed

\vskip 10pt
Families ${\cal F}$ for which $\forb(m,{\cal F})$ is $O(1)$ must arise from applying \lrf{FinG} and  \trf{constant} in view of \trf{classify}.
There are no 2-fold or 1-fold product constructions in common for $Q_1,Q_2$ so that $X(\{Q_1,Q_2\})=1$. We can use \trf{constant} and \lrf{FinG} to get a constant bound but perhaps recording a general result is in order. Let $0_{a,b}$ denote the $a\times b$ matrix of 0's and let $J_{a,b}$ denote the $a\times b$ matrix of 1's.

\begin{thm}  Let $k, \l, p, q$ be given.  Then there exists some constant $c_{k\l pq}$ such that for $m \geq c_{k\l pq}$, we have $\forb(m, \{0_{k,\l}, J_{p,q}\}) = \l + q - 2$. \label{thm:constant}\end{thm}

\proof We let $d = \max\{k,\l,p,q\}$.  Then  $0_{k,\l} \prec T_{2d}, 0_{k,\l} \prec I_{2d}$ and $ J_{p,q} \prec I^c_{2d}$.  Thus by \trf{classify}, 
$\forb(m, \{0_{k,\l}, J_{p,q}\})$ is $O(1)$.
 We wish to show that $\forb(m, \{0_{k,\l}, J_{p,q}\}) = \l + q - 2$.  Let $B\in\Av(m, \{0_{k,\l}, J_{p,q}\})$ with $n=\ncols{B}>\l+q-2$. We can delete columns if necessary to obtain a matrix $A\in\Av(m, \{0_{k,\l}, J_{p,q}\})$ with $n=\ncols{A}=\l+q-1$.
  From \lrf{lem:sum} we know that the right side of \rf{eqn:sum} is constant based on $n, k, \l, p$ and $q$. The right hand side of the inequality in \rf{eqn:sum} is at least $m$ since the summands of the left side will be at least 1 unless  $a_r < \l$ and $b_r < q$ which is impossible  because $a_r + b_r = \l+q-1$.   So for sufficiently large $m$, we have a contradiction. Hence there exists a constant $c_{k\ell pq}$ so that for $m \geq c_{k\l pq}$, we have $\forb(m, \{0_{k,\l}, J_{p,q}\}) \le \l + q - 2$.

It remains to show we have a construction $A\in\Av(m,\{0_{k,\l}, J_{p,q}\})$ with $\ncols{A}=\ell+q-2$.   Assume $m=\binom{\ell+q-2}{q-1}+t$ for some $t\ge0$. Let the first $\binom{\ell+q-2}{q-1}$ rows of $A$ consist of all possible rows of $\ell+q-1$ entries with exactly $q-1$ 1's. For the remaining rows of $A$ simply repeat the row of $q-1$ 1's followed by $\ell-1$ 0's $m-\binom{\ell+q-2}{q-1}$ times. The matrix is seen to be simple and cannot have $0_{k,\ell}$ since each row has $\ell-1$ 0's and cannot have $J_{p,q}$ since each row has $q-1$ 1's. Thus  $\forb(m, \{0_{k,\l}, J_{p,q}\}) \ge q + \l - 2$.  This yields the result.
\qed

\begin{lemma} Let $k, \l, p, q$ be given.  Let $A \in \Av(m, \{0_{k,\l}, J_{p,q}\})$, with $\ncols{A} = n$.  Also let $a_r$ denote the number of 0's in row $r$ of $A$, and $b_r$ the number of 1's in row $r$ so that $a_r + b_r = n$.  Then:
\begin{equation}\sum_{r=1}^m \left(\binom{a_r}{\l}+\binom{b_r}{q}\right) \leq (k-1)\binom{n}{\l} + (p-1)\binom{n}{q}.  \label{eqn:sum}\end{equation}
\label{lem:sum}\end{lemma}

\proof We consider the columns of $A$.  We take all $\l$-subsets of the columns and call them $0$-buckets.  Similarly, we take all $q$-subsets of the columns as $1$-buckets.  We will have $\binom{n}{\l}$ $0$-buckets and $\binom{n}{q}$ $1$-buckets.  We then process the rows of $A$ one by one, considering all possible $\l$-subsets and $q$-subsets of columns on that row.  If one of these subsets contains all 0's or all 1's, it makes a contribution to the appropriate $0$-bucket or $1$-bucket.  Thus if there are $a$ 0's in a row, and $b$ 1's (where $a + b = n$), then the row will make contributions to $\binom{a}{\l}$ $0$-buckets and $\binom{b}{q}$ $1$-buckets.  The left side of \rf{eqn:sum} is thus the total number of contributions over the rows of $A$.  Each of our $\binom{n}{\l}$ $0$-buckets can have a maximum of $k-1$ contributions, and similarly, our $\binom{n}{q}$ $1$-buckets can have a maximum of $p-1$ contributions, which produces the right side of the inequality.  
\qed 
\vskip 5pt

\begin{thm} We have that $\forb(m,\{Q_1,Q_2\})$, $\forb(m,\{Q_1,Q_5\})$, $\forb(m,\{Q_1,Q_7\})$,\hfil\break $\forb(m,\{Q_2,Q_4\})$, $\forb(m,\{Q_2,Q_6\})$, $\forb(m,\{Q_4,Q_5\})$, $\forb(m,\{Q_4,Q_7\})$ and  \hfil\break $\forb(m,\{Q_5,Q_6\})$
 are all bounded by $O(1)$.\label{constantbd}\end{thm}

\proof  We apply \lrf{FinG} with ${\cal G}=\{I_4,I_4^c,T_4\}$ and also \trf{constant}. Two examples are the following. For the family $\{Q_1,Q_5\}$ we note that $Q_1\prec I_4$, $Q_5\prec I_4^c$ and $Q_1\prec T_4$. For the family $\{Q_5,Q_6\}$ we note that $Q_6\prec I_4$, $Q_5\prec I_4^c$ and $Q_5\prec T_4$.\qed

\vskip 10pt
We record below the exact values for $\forb(m,\{Q_1,Q_2\})$. The function $\forb(m,\{Q_1,Q_2\})$  has a surprising non-monotonicity in $m$.
\begin{thm} \cite{thesis} We have
$$
\forb \left(m, \left\{ \left[
              \begin{array}{cc}
                0 & 0 \\
                0 & 0 \\
              \end{array}
            \right], \left[
              \begin{array}{cc}
                1 & 1 \\
                1 & 1 \\
              \end{array}
            \right]
 \right\} \right) = \begin{cases}
  2 & \text{if } m=1 \hbox{ or } m\geq7 \\
  4  & \text{if } m = 2,5,6 \\
  6 & \text{if } m = 3,4
  \end{cases}.
$$
\label{thm:failedmono}\end{thm}

\section{Graph Theory}\label{graphtheory}

We consider a family ${\cal F}=\{\1_3, F\}$ for some $F$. Note that $Q_5=\1_3$.  We have that $\forb(m,\{\1_3,F\})$ is $O(m^2)$ since $\forb(m,\1_3)$ is $O(m^2)$. In this section we consider those $F$ which are (0,1)-matrices with column sums 0,1 or 2.  If $F$ has a repeated column of sum 2 then $2\cdot \1_2\prec F$ and then $\forb(m,\{\1_3,F\})$ is $\Theta(m^2)$ (the construction $I_{m/2}\times I_{m/2}$ yields the lower bound).  So we may assume $F$ has no repeated columns of sum 2 and so these columns can be viewed as the incidence matrix of some graph.   We will adapt \lrf{extendgraph} to those $F$ with columns having sum 0,1 or 2.  The following remark describes our construction.

\begin{remark} Let $F$ be a $k\times \ell$  (0,1)-matrix with column sums $\in \{0,1,2\}$. Assume $2\cdot \1_2\not\prec F$. Let $a_i$ be the number of columns of $F$ of sum 1 with a 1 in row $i$, and let $b$ be the number of columns of $F$ of all 0's.  We can form a graph $G$ with $V(G)=[k+\sum_{i\in [k]}a_i+b+1]$ as follows. For $i,j\in [k]$ we have $i,j\in E(G)$ if and only if there is a column of $F$ with 1's in rows $i,j$. Also, for each $i\in [k]$, we add $a_i$ edges to $G$ joining $i\in [k]$ to $a_i$ vertices chosen from
$[k+\sum_{i\in [k]}a_i+b+1]\backslash [k]$ (each of which has degree 1). Finally on the remaining $b+1$ vertices we add $b$ edges in the form of a tree. Then $F\prec \inc(G)$. \label{adorn}\end{remark}
\proof  We find $F\prec \inc(G)|_{[k]}$. \qed

\vskip 10pt
The remark demonstrates some of the differences between a `subgraph' and a `configuration'.

\begin{lemma} Let $T$ be a graph on $k$ vertices and assume $T$ has  no cycles (i.e. a forest). Then $\ex(m,T)$ is $O(m)$. \label{treebd}\end{lemma}
\proof  Folklore  says if a graph $G$ on $m$ vertices has at least $km$ edges then $T$ is a subgraph of $G$.  Assume $G$ has at least $km$ edges. We first obtain a subgraph $G'$ of $G$ with minimum degree $k$ which we obtain by removing  vertices whose degree is at most $k-1$. Each vertex deleted removes at most $k-1$ edges. Thus the process must stop with a non-empty subgraph $G'$ of $G$ with minimum degree $k$.  Since $T$ has no cycles, we may order the vertices $v_1,v_2,\ldots, v_k$ of $T$ so that for each $v_i$ there is at most one $v_j$ with $j<i$ such that $(v_j,v_i)$ is an edge of $T$. Assume we have found in $G'$ a subgraph $T$, namely vertices $x_1,x_2,\ldots ,x_p\in V(G')$ such that  $(x_i,x_j)\in E(G')$  if $(v_i,v_j)\in E(T)$ where $1\le i<j\le p$. If $p=k$, we are done. If $p<k$, then consider $v_{p+1}$. If $v_{p+1}$ is not joined in $T$ to anything in $v_1,v_2,\ldots ,v_p$ then we can select $x_{p+1}$ as any  vertex in $G'$ (say adjacent to $x_p$) which has not already been selected. 
 If $v_{p+1}$ is joined to $v_i$ with $i\le p$, then we choose $x_{p+1}$ as any vertex adjacent $x_i$ which has not already been selected. We use that minimum degree in $G'$ is at least $k>p$. Continue until $p=k$.  We deduce that $\ex(m,T)<km$. \qed
 
 \vskip 10pt
 We  extend this to configurations in \trf{graphlinear}. 
 \begin{thm} Let $F$ be  a given $k\times \ell$ (0,1)-matrix such that every column has at most 2 1's. Assume that $2\cdot\1_2\not\prec F$ and assume $C_t\not\prec F$ for any $t\ge 3$.  
 Then $\forb(m,\{\1_3,F\})$ is $O(m)$.\label{graphlinear}\end{thm}
 
 \proof Use \remrf{adorn} to obtain a graph $G$ from $F$.  We check that $G$ has no cycles and hence by our above remarks, $\ex(m,G)$ is $O(m)$.  We note that $F\prec \inc (G)$.  Now applying  \lrf{extendgraph} yields 
 $\forb(m,\{\1_3,\inc(G)\})$ is $O(m)$ and so, by \lrf{FinG},   $\forb(m,\{\1_3,F\})$ is $O(m)$.\qed
 
 \vskip 10pt
The following is a weak version of the extremal graph results of Erd\H{o}s, Stone and Simonovits since we only consider asymptotic growth rates.

\begin{thm} Let $F$ be  a given $k\times \ell$ (0,1)-matrix such that every column has at most 2 1's. Let $t$ be given. Assume $2\cdot\1_2\prec F$ or there is some $t\ge 1$ with $C_{2t+1}\prec F$.  Then
$\forb(m,\{\1_3,F\})$ is $\Theta(m^2)$.\label{oddcycle}\end{thm}

\proof The upper bound $O(m^2)$ is easy. We may use the construction $I_{m/2}\times I_{m/2}$ to obtain the matching lower bound.\qed

\vskip 5pt
Let $H$ be a bipartite graph. Then $\ex(m,H)$ is $o(m^2)$. We use the notation $o(m^2)$ to refer to a function $f(m)$ with
$\lim_{m\to\infty}f(m)/m^2=0$.  The following result extends this to configurations. 

\begin{thm} Let $F$ be  a given $k\times \ell$ (0,1)-matrix such that every column has at most 2 1's.
Let $F$ be given with  and also with the property that  $2\cdot\1_2\not\prec F$ and for all $t\ge 1$, we have $C_{2t+1}\not\prec F$.  Then
$\forb(m,\{\1_3,F\})$ is $o(m^2)$.\end{thm}

\proof Form a graph $G$ as described in \remrf{adorn}. Since for all $t\ge 1$, we have $C_{2t+1}\not\prec F$, the resulting graph $G$ will be a bipartite graph.  Then for some $s,t$, $G$ is a subgraph of the complete bipartite graph $K_{s,t}$. We know that
$\ex(m,K_{s,t})$ is $o(m^2)$. Thus $\ex(m,G)$ is  $o(m^2)$. Now $F\prec \inc(G)$ and so by \lrf{extendgraph} we have that  $\forb(m,\{\1_3,F\})$ is $o(m^2)$.\qed

\vskip 5pt
One could imagine trying to obtain similar results for $\forb(m,\{\1_k,F\})$ where $F$ has columns sums at most $k-1$. It is still very much an open problem to determine the exact asymptotic growth $\ex(m,C_{2t})$ for various $t\ge 2$ with two results noted \trf{C4}, \trf{C6}.   \trf{graphlinear} combined with \remrf{complement}, yields the following.

\begin{thm} We have that $\forb(m,\{Q_5,Q_3\})$,  $\forb(m,\{Q_5,Q_8\})$,  $\forb(m,\{Q_5,Q_9\})$,  \hfil\break 
 $\forb(m,\{Q_4,Q_3\})$, $\forb(m,\{Q_4,Q_8\})$, 
$\forb(m,\{Q_4,Q_9\})$ are all $O(m)$.\label{graphbd}\end{thm}

\trf{oddcycle} yields that $\forb(m,\{Q_5,Q_7\})$ is $\Omega(m^2)$, a fact which has already been noted.

\section{New Standard Induction}\label{SI}
Our standard induction argument proceeds as follows. Let  $A\in\Av(m,{\cal F})$ with $\ncols{A}=\forb(m,{\cal F})$. We choose  $r\in [m]$ and delete row $r$ from $A$. The result may have repeated columns which we collect in a matrix $C_r$. After permuting rows and columns we have the following:
\begin{E}
A = 
\begin{array}{r@{}}
    \hbox{row  }r \\  \\
\end{array}
\begin{bmatrix}
    0 \,0 \cdots  0 & 1\,1  \cdots  1\\
    B_r  \,\, C_r  & C_r  \,\, D_r \\
\end{bmatrix}.
\label{rdecomp}\end{E}
Both $[B_r\,C_r\,D_r]$ and $C_r$ are simple. We have $[B_r\,C_r\,D_r]\in \Av(m-1,{\cal F})$ suggesting an induction.  
Now $[0\,1]\times C_r$ is in $A$. Thus define ${\cal G}$ as the \emph{minimal} set of configurations $F'$ such that $F\prec [0\,1]\times F'$ for some $F\in{\cal F}$ (we defined minimal after \remrf{minimal}). We deduce that $C_r\in\Av(m-1,{\cal G})$.
This yields the following induction formula
 \begin{E}\forb(m,F,s)= \ncols{A}=\ncols{[B_rC_rD_r]}+\ncols{C_r}\leq \forb(m-1,{\cal F})+ \forb(m-1,{\cal G}).\label{recursion}\end{E}
This means any upper bound on $\ncols{C_r }$ (as a function of $m$) automatically yields an upper bound on $A$ by induction. Thus to show $\forb(m,{\cal F})$ is $O(m)$ it suffices to show $\ncols{C_r}$ is bounded by a constant.  
We have discovered a new standard induction \cite{AL} that, by extending the argument to matrices with multiple columns, yields a more powerful induction formula \rf{newrecursion}.   
Let $A$ be an m-rowed (0,1)-matrix (not necessarily simple) and $\alpha$ be an $m\times 1$ column. Let $\mu(\alpha,A)$ denote the multiplicity of column $\alpha$ in $A$.
We say $A$ is $s$-\emph{simple} if every column $\alpha$ of $A$ has $\mu(\alpha,A)\le s$. Let $\Av(m,{\cal F},s)$ denote the $m$-rowed $s$-simple matrices with no 
$F\in{\cal F}$. We define 
$$\forb(m,{\cal F},s)=\min_A\{\ncols{A}\,:\,A\in\Av(m,{\cal F},s)\}.$$
We note that $\forb(m,{\cal F})\le\forb(m,{\cal F},s)\le s\cdot\forb(m,{\cal F})$ and so the asymptotic growth rate of $\forb(m,{\cal F})$ and $\forb(m,{\cal F},s)$ are the same (for fixed $s$).  Associate with $A$ the simple matrix $\supp(A)$ where $\mu(\alpha,\supp(A))=1$ if and only if $\mu(\alpha,A)\ge 1$.
Given ${\cal F}$, let $t$ be the maximum multiplicity of a column in $F$ over all $F\in{\cal F}$, i.e. each $F\in{\cal F}$ is $t$-simple but some $F\in{\cal F}$ is not $(t-1)$-simple.  We assume for \rf{newrecursion} that some $F\in{\cal F}$ is not simple and so $t\ge 2$.  Define $s=t-1$.  Assume $A\in\Av(m,{\cal F},s)$. We first decompose $A$ using row $r$ as follows:
$$A = 
\begin{array}{r@{}}
    \hbox{row  }r \\  \\
\end{array}
\begin{bmatrix}
    0 \, 0\cdots  0 & 1 \,1 \cdots  1\\
    G&H \\
\end{bmatrix}.
$$
We deduce that $\mu(\alpha,G)\le s$ and $\mu(\alpha,H)\le s$.  We can obtain the following decomposition of $A\in\Av(m,{\cal F},s)$ based on deleting row $r$ and rearranging by selecting certain columns for $C_r$ so that   if $\mu(\alpha,G)+\mu(\alpha,H)\ge s+1$, then  $\mu(\alpha,C_r)=\min \{\mu(\alpha,G),\mu(\alpha,H)\}$.  We again obtain \rf{rdecomp}.
We conclude that $[B_r\,C_r\, D_r]$ and $C_r$ are both $s$-simple. Thus $\ncols{[B_r\,C_r\, D_r]}\le\forb(m,{\cal F},s)$.
 Since each column in $C_r$ appears at least $s+1$ times in $[B_rC_rC_rD_r]$, then $C_r$  has no configuration in ${\cal F}'=\{\supp(F)\,:\,F\in{\cal F}\}$.
 In the case that each  $F\in{\cal F}$ is  simple then  ${\cal F}'={\cal F}$.  
 We obtain the following useful inductive formula:
 \begin{E}\forb(m,F,s)= \ncols{A}=\ncols{[B_rC_rD_r]}+\ncols{C_r}\leq s\cdot\left(\forb(m-1,{\cal F})+ \forb(m-1,{\cal F}'\cup{\cal G})\right).\label{newrecursion}\end{E}
 The extra value here as compared with \rf{recursion} is in forbidding in $C_r$ the configurations $\supp(F)$ for each $F\in {\cal F}$.

\begin{thm}\label{F8blockof0's}Let $k,\ell$ be given. Then $\forb(m,\{Q_8,[0\,1]\times 0_{k,\ell}\})$ is $O(m)$.\end{thm}
\proof Let $A\in\Av(m,\{Q_8,[0\,1]\times 0_{k,\ell}\})$.  We apply the decomposition of \rf{rdecomp} and deduce that
$C_r\in\Av(m-1,\{I_2,0_{k,\ell}\})$. We note that $Q_8=[0\,1]\times I_2$ and deduce that ${\cal G}=\{I_2,0_{k,\ell}\}$. With $I_2\not\prec C_r$, we discover that 
$C_r\prec [\0_{m-1}|T_{m-1}]$ (i.e. $C_r$ is a selection of columns from the triangular matrix).  Then if $\ncols{C_r}\ge k+\ell$, we find  $0_{k,\ell}\prec C_r$. We deduce that $\ncols{C_r}\le k+\ell-1$ and  deduce by induction on $m$ (using \rf{recursion}) that 
$\forb(m,\{Q_8,[0\,1]\times 0_{k,\ell}\})$ is $O(m)$.\qed

\begin{cor}We have that $\forb(m,\{Q_1,Q_8\})$, $\forb(m,\{Q_2,Q_8\})$, $\forb(m,\{Q_4,Q_8\})$ and $\forb(m,\{Q_5,Q_8\})$  are  $O(m)$.\label{F1F8,F2F8}\end{cor}

\proof We note that $Q_1\prec [0\,1]\times 0_{1,2}$ and $Q_4\prec [0\,1]\times 0_{2,1}$.  Also $Q_8^c$ is the same configuration as $Q_8$ and $Q_1^c=Q_2$, $Q_4^c=Q_5$ so we apply \remrf{complement}.  \qed

\begin{thm}Let $t\ge 2$ be given. Then $\forb(m,\{Q_8, t\cdot([0\,1]\times [0\,1])\})$ is $O(m)$.\end{thm}
\proof Let $A\in\Av(m,\{Q_8,t\cdot([0\,1]\times [0\,1])\})$.  We apply the decomposition obtained as \rf{rdecomp} and deduce that
$C_r\in\Av(m-1,\{I_2,t\cdot [0\,1]\}$. We note that $Q_8=[0\,1]\times I_2$ and deduce that ${\cal G}=\{I_2,t\cdot [0\,1]\}$.   With $I_2\not\prec C_r$, we have that $C_r\prec [\0_{m-1}|T_{m-1}]$. 
For $\ncols{C_r}\ge 2t$, we deduce that $t\cdot [0\,1]\prec C_r$. This is a contradiction and so $\ncols{C_r}\le 2t-1$ and  deduce by induction on $m$ (using \rf{recursion}) that   $\forb(m,\{Q_8, t\cdot([0\,1]\times [0\,1])\})$ is $O(m)$.
\qed

\vskip 5pt
We note that $Q_3\prec 2\cdot([0\,1]\times [0\,1])$ and obtain the following. 

\begin{cor}We have that $\forb(m,\{Q_3,Q_8\})$ is $O(m)$.\label{Q3Q8}\end{cor}

We have that $Q_3\not\prec I\times I^c$ and $Q_3$ is a configuration in the other five 2-fold products.  
We have that $Q_6\not\prec I^c\times I^c$, $Q_6\not\prec I^c\times T$ and $Q_6\not\prec T\times T$. Also $Q_6$ is a configuration in the other three 2-fold products. 	 We also note that either $T$, $I^c$ are 1-fold products avoiding $Q_3$ and $Q_6$. 
Let 
$$F_2(1,t,t,1)=\left[\begin{array}{c}0\\ 0\\ \end{array}\right.
\overbrace{  
\begin{array}{c}1\,1\cdots 1\\ 0\,0\cdots 0\\ \end{array}
}^t
\overbrace{  
\begin{array}{c}0\,0\cdots 0\\ 1\,1\cdots 1\\ \end{array}
}^t
\left.\begin{array}{c}1\\ 1\\ \end{array}\right].$$
We are using notation from \cite{survey}. Thus $Q_3=F_2(1,2,2,1)$. We have that $F_2(1,t,t,1)\not\prec I\times I^c$ and  $F_2(1,t,t,1)$ is a configuration in the other five 2-fold products.   Similarly to $Q_6=I_3$, the configuration $t\cdot I_k$ is not in the $(k-1)$-fold products consisting solely of the terms  $I^c$ and $T$ but is in every 2-fold product using $I$.
Thus we might guess (using \corf{grand})  that forbidding $F_2(1,t,t,1)$ and $t\cdot I_k$  results in a linear bound. This is true.
  The following  is proven using two lemmas.

  \begin{thm}\label{tIkF2(1,t,t,1)} Let $k,t\ge 2$ be given.
 We have that $\forb(m,\{t\cdot I_k,F_2(1,t,t,1)\})$ is $\Theta(m)$. \end{thm}
 
 \proof  We will use induction on $m$. Let $A\in\Av(m,\{t\cdot I_k,F_2(1,t,t,1)\},t-1)$. We use $s=t-1$ and then for any row $r\in [m]$, obtain the decomposition \rf{rdecomp}.
 We wish to use \rf{newrecursion}.  With ${\cal F}=\{t\cdot I_k,F_2(1,t,t,1)\}$ and $s=t-1$, we have ${\cal F}'=\{I_k,F_2(1,1,1,1)\}$ (since $F_2(1,t,t,1)\prec t\cdot F_2(1,1,1,1)$) and ${\cal G}=\{t\cdot [1\,0], t\cdot [\0_{k-1}\,|\,I_{k-1}]\}$ (since $F_2(1,t,t,1)\prec [0\,1]\times (t\cdot [0\,1])$).  
 Then $C_r\in\Av(m,\{I_k,F_2(1,1,1,1),t\cdot [1\,0], t\cdot [\0_{k-1}\,|\,I_{k-1}]\},t-1)$.   The second configuration ($F_2(1,1,1,1)$) and the fourth configuration $( t\cdot [\0_{k-1}\,|\,I_{k-1}])$ do not get used in our proof.
 We form a digraph on $[m]$ by setting $r\rightarrow s$ if there are at most $t-1$ columns of $A$ with $\linelessfrac{r}{s}\left[\linelessfrac{1}{0}\right]$. If there is a row $s$ of $C_r$ with one 0 and at least $t$ 1's then, by considering the forbidden configuration $F_2(1,t,t,1)$,  we deduce that $r\rightarrow s$ (else $F_2(1,t,t,1)\prec A|_{\{r,s\}}$). Given a row $r$, assume no such row $s$ exists.
 Then all rows of $C_r$ have either at most $t-1$ 1's or is all 1's.
 
 Assume $\ncols{C_r}\ge tk$. Now remove from $C_r$ any rows of all 1's to obtain a simple matrix $C'$ and obtain a simple matrix $C$ from $C'$ by deleting a column of 0's if it exists. We deduce that each row of $C$ has  at most $t-1$ 1's and each column of $C$ has at least one 1. Also $\ncols{C}\ge tk-1\ge (t-1)k$.  By \lrf{findIk}, we deduce that $C_r$ has $I_k$, a contradiction.  So a row $s$ exists. Since for each row $r\in [m]$ there is some row $s\in [m]$ with 
 $r\rightarrow s$, we deduce that there is a directed cycle and we may  apply \lrf{cycleoffalls}  to show that $\ncols{A}$ is $O(m)$.\qed
 
 \vskip 10pt
 We have used the following idea before.
\begin{lemma}\label{cycleoffalls} Let $A$ be a simple matrix for which there are $k$ rows $a_1,a_2,\ldots ,a_k$ for which there are at most $t$ columns containing 
$\linelessfrac{a_i}{a_{i+1}} \bigl[\linelessfrac{1}{0}\bigr]$ for $i=1,2,\ldots ,k-1$ and also   there are at most $t$ columns containing  
$\linelessfrac{a_k}{a_{1}} \bigl[\linelessfrac{1}{0}\bigr]$. 
Then we may delete  up to $kt$ columns from $A$ (as described) and the $k-1$ rows $a_1,a_2,\ldots ,a_{k-1}$  and obtain a simple matrix. \end{lemma}

\proof 
Consider the matrix $A'$ obtained from $A$ by deleting the special columns described of which there are at most  $kt$. Then we deduce that $A'|_{\{a_1,a_2,\ldots ,a_k\}}$ consists of columns of all 0's and columns of all 1's.  Now deleting from $A$ the special columns and the  $k-1$ rows  $a_1,a_2,\ldots ,a_{k-1}$ will result in a simple matrix. \qed

\begin{lemma} Let $C$ be a matrix having row sums at most $t-1$. Assume each column sum is at least 1. Assume $\ncols{C}\ge (t-1)k$, Then $I_k\prec C$.
\label{findIk}\end{lemma}
\proof We could phrase this with sets corresponding to the rows.  For row $r$ we form a subset $S_r\subseteq \{1,2,\ldots ,\ncols{C}\}$ with $s\in S_r$ if and only if there is a 1 in row $r$ and column $s$. Our induction is on $k$ with the result being trivial for $k=1$.  We can greedily select sets  $S_1,S_2,\ldots S_p$ so that 
$S_j\backslash \left(\cup_{i=1}^{j-1}S_i\right)\ne\emptyset$ for $j\in[p]$ and 
so that $\cup_{i=1}^pS_i\ge (t-1)k$.   We begin by choosing an element $a_k\in S_p\backslash \left(\cup_{i=1}^{p-1}S_i\right)\ne\emptyset$. We delete, from our  $p$ sets, the 
 elements of $S_p$ (there are at most $t-1$ such elements) and then delete any sets which are now $\emptyset$. We now have sets $S_1',S_2',\ldots S_q'$ so that $S_i'\backslash \cup_{i=1}^{j-1}S_i'\ne\emptyset$ and 
so that $\cup_{i=1}^qS_i\ge (t-1)(k-1)$ and $|S_i|\le t-1$. If we form a set-element incidence matrix $C'$ from these $q$ sets, we find that each row sum of $C'$
 is at most $t-1$ ($|S_i'|\le |S_i|\le t-1$). Moreover 
each column sum is at least 1 (we deleted columns corresponding to elements of $S_p$) and $\ncols{C}\ge (t-1)(k-1)$  (we only deleted the elements of $S_p$ and $|S_p|\le t-1$). By induction on $k$, $ I_{k-1}\prec C'$. Now the $p$th row of $C$ is 0's on columns not in $S_p$ and in column $a_k$ has 0's on all rows except row $p$ for which it is 1.  Now we find $I_k\prec C$.
\qed
\vskip 10pt
\begin{cor}We have that $\forb(m,\{Q_1,Q_3\})$, $\forb(m,\{Q_2,Q_3\})$, $\forb(m,\{Q_3,Q_6\})$ and $\forb(m,\{Q_3,Q_7\})$ are $O(m)$.\label{orderbd}\end{cor}
\proof We use \lrf{FinG} with ${\cal G}=\{F(1,t,t,1),t\cdot I_k\}$. For example $Q_1\prec t\cdot I_k$ and $Q_3\prec F(1,t,t,1)$ and also
$Q_6\prec t\cdot I_k$.  We also use \remrf{complement} noting that $\{Q_1^c,Q_3^c\}$ and $\{Q_2,Q_3\}$ are the same as sets of configurations and 
$\{Q_3^c,Q_6^c\}$ and $\{Q_3,Q_7\}$ are the same as sets of configurations. \qed

\section{Structure that arises from forbidding $Q_9$}\label{structure}

The following result gives some of the structure of matrices $A\in\Av(m,Q_9)$.
Let $A_k$ denote the  columns of  $A$ of column sum $k$.  We discover that $A_k$ is of one of two types.
 We say $A_k$ is of \emph{type} 1 if there is a partition of the rows
$[m]=X_k\cup Y_k\cup Z_k$ such that all columns in $A_k$ are 1's on rows $X_k$, 0's on rows $Z_k$ and each column of $A_k$ has exactly one 1 in rows $Y_k$
Thus $A_k|_{Y_k}$ is $I_{|Y(k)|}$. In that case, by examining column sums, $|X_k|+1=k$. We say $A_k$ is of \emph{type} 2 if there is a partition of the rows
$[m]=X_k\cup Y_k\cup Z_k$ such that all columns in $A_k$ are 1's on rows $X_k$, 0's on rows $Z_k$ and each column of $A$ has exactly one 0 in rows $Y_k$
Thus $A_k|_{Y_k}$ is $I_{|Y(k)|}^c$. In that case, by examining column sums, $|X_k|+|Y_k|-1=k$.

\begin{lemma}\cite{FFP87} Let $A\in\Av(m,Q_9)$. Let $A_k$ denote the  columns of column sum $k$. Then $A_k$ is of type 1 or type 2.\label{types}\end{lemma}
\proof  This follows quite readily by considering the columns of $A_k$ one column at a time. For $\ncols{A_k}\le 2$, the type would not be unique. \qed

\vskip 10pt
We consider the following $(t+1)\times (2t+2)$ matrix $F(t)$ whose first two rows coincide with $F_2(1,t,t,1)$:
$$F(t)=\left[\begin{array}{@{}c@{}}\\ \\ \\ \\ \end{array}\right. 
\begin{array}{@{}c}0\\0\\ \vdots\\ 0\\ \end{array}
\overbrace{\begin{array}{c}1\,1\cdots 1\\0\,0\cdots 0\\ \vdots\\ 0\,0\cdots 0\\ \end{array}}^t
\overbrace{\begin{array}{c}0\,0\cdots 0\\1\,1\cdots 1\\ \vdots\\ 1\,1\cdots 1\\ \end{array}}^t
\left.\begin{array}{c}1\\1\\ \vdots\\ 1\\ \end{array}\right].$$

\begin{lemma} Let $t\ge 1$ be given. Then  $\forb(m,\{Q_9,F(t)\})$ is $O(m)$.\label{F(t)}\end{lemma}
\proof Let $A\in\Av(m,\{Q_9,F(t)\})$. We will show that $\ncols{A}\le (7t+1)m$.
let $A_k$ denote the  columns of column sum $k$.  For $j=1,2$, let $W(j)=\{k\,:\,A_k\hbox{ is of type }j, \ncols{A_k}\ge t+2\}$ and let $V(j)$ be the concatenation of $A_{k}$ for $k\in W(j)$ so that $\ncols{V(j)}=\sum_{k\in W(j)}\ncols{A_k}$. 
 
We first note that for $a<b$ that $|X_a\backslash X_b|\le 1$. This is because if $|X_a\backslash X_b|\ge 2$ and $r,s\in X_a\backslash X_b$ then any column $\alpha$ from $A_a$ has  1's on rows $r,s$. We can choose a column $\beta$ from $A_b$ with  0's  on rows $r,s$ using $r,s\in Y_b\cup Z_b$ and the fact that $\ncols{A_b}\ge t+2$.  But $\beta$ has more 1's than $\alpha$ and so we find $Q_9\prec [\alpha\,|\,\beta]$. 

Assume $\ncols{V(1)}\ge 3tm+1$. Then there are $3t$ indices $\{s(1),s(2),\ldots ,s(3t)\}\subseteq W(1)$ where $s(1)<s(2)<\cdots <s(3t)$ so that we can find row $r$ with  $r\in\cap_{i=1}^{3t}Y_{s(i)}$.  We wish to find a set of rows $S$ with $|S|=t$ such that
$S\subseteq X_{s(3t)}\cap\left(\cup_{i=1}^t(Y_{s(i)}\cup Z_{s(i)})\right)$.  We have $|X_{s(3t)}\backslash X_{s(t)}| \ge 2t$.  Using $|X_{s(i)}\backslash X_{s(t)}|\le 1$, we have
$|X_{s(3t)}\backslash \left(\cup_{i=1}^tX_{s(i)}\right)| \ge t$ and so we can find $S$ as claimed.
Now we obtain $F(t)$ as follows. For each $i$ with $1\le i\le t$, we have $r\in Y_{s(i)}$ and $S\subseteq Y_{s(i)}\cup Z_{s(i)}$. We choose one column from $A_{s(1)}$ with a 0 on row $r$ where we choose the column so it also has 0's on rows $S$ (which is possible for  $\ncols{A_{s(i)}}\ge t+2$ (else with $r\cup S=Y_{s(i)}$ we would have difficulty finding the column). We choose one column from each $A_{s(i)}$ for $i\in[t]$, with a 1 on row $r$ and necessarily 0's on rows $S$ and one    All columns from $A_{s(3t)}$ are 1's on rows $S\subseteq X_{s(3t)}$. With $\ncols{A_{s(i)}}\ge t+2$, we can find $t+1$ columns in $A_{s(3t)}$ of which $t$ are 0 on row $r$ and one is 1 on row $r$.   This completes $F(t)$. We conclude that $\ncols{V(1)}\le 3tm$
 Noting that $Q_9^c,F(t)^c$ are the same as $Q_9,F(t)$ when considered as configurations, we deduce that  $\ncols{V(2)}\le 3tm$.
Now $A$ consists of $V(1)$ and $V(2)$ plus at most $(t+1)m$ columns (to account for $\ncols{A_k}$ where $\ncols{A_k}\le t+1$) and so $\ncols{A}\le  (7t+1)m$. \qed

\begin{cor}We have that $\forb(m,\{Q_1,Q_9\})$,  $\forb(m,\{Q_2,Q_9\})$, $\forb(m,\{Q_3,Q_9\})$,\break $\forb(m,\{Q_4,Q_9\})$ and $\forb(m,\{Q_5,Q_9\})$ are $O(m)$.\label{Q9}\end{cor}
\proof  We note that $Q_1, Q_2, Q_3, Q_4, Q_5$ are all configurations in $F(2)$. 
\qed

\begin{thm}We have that $\forb(m,\{Q_8,Q_9\})$ is $O(m)$.\label{Q8Q9}\end{thm}
\proof  Let $A\in\Av(m,\{Q_8,Q_9\})$ and let $A_k$ denote the  columns of column sum $k$. Assume $\ncols{A_k}\ge 3$ for all $k$. 
For $j=1,2$, let $W(j,\even)=\{k\,:\,A_k\hbox{ is of type }j, \ncols{A_k}\ge 3, j\hbox{ is }\even\}$ and let $V(j, \even)$ be the concatenation of $A_{k}$ for $k\in W(j,\even)$.  
We similarly define $W(j,\odd)$ and $V(j,\odd)$. This more complicated definition ensures that for $a,b\in W(j,\even)$ (or $a,b\in W(j,\odd)$) with $a<b$  that $a<a+1<b$ (column sums differ by at least 2).

We wish to show $|V(1,\even)|\le 2m$. We establish a number of properties before using an interesting induction.
  We may assume that for $i<j$ and $i,j\in W(1,\even)$, that $|X_i\backslash X_j|\le 1$ else we have a copy of $Q_9$ in $[A_i\,|\,A_j]$ as described in proof of \lrf{F(t)}.

We may assume $|Y_i\cap Y_j|\le 1$ for all pairs $i,j\in W(1,\even)$.  Otherwise assume $|Y_i\cap Y_j|\ge 2$ for some pair $i<j$ with $i,j\in W(1,\even)$. Let $r,s\in Y_i\cap Y_j$. Now $|X_i|<|X_j|$ and so we can choose a third row $p\in X_j\backslash X_i$.
We now find a copy of $Q_8$ in $[A_i\,|\,A_j]$ in rows $p,r,s$, a contradiction

Now assume $|Y_i\cap Y_j|= 1$ for some pair $i<j$.   We claim $X_i\subset X_j$. Otherwise, choose $r\in X_i\backslash X_j$ and $p=Y_i\cap Y_j$. We can find
$\left[\linelessfrac{1}{1}\right]$ in rows $p,r$ of some column of $A_i$ and  $\left[\linelessfrac{0}{0}\right]$ in rows $p,r$ of some column of $A_j$. Give $i<j$ we now have a copy of $Q_9$ in $[A_i\,|\,A_j]$, a contradiction.  

Finally assume we have indices 
$a,b,c\in W(1)$ with $a<b<c$ and $Y_a\cap Y_c=\{r\}$ and $Y_b\cap Y_c=\{s\}$. Then we conclude $r=s$. If not,
 recall that $X_a\subset X_c$ and $X_b\subset X_c$ and $a<a+1<b<b+1<c$. Now $|X_c\backslash X_b|\ge 2$ and  $|X_{a}\backslash X_b|\le 1$ so we are able choose $p\in X_{c}\backslash (X_{b}\cup X_a)$. Then we find $Q_8$ in rows $p,r,s$ of  $[A_a\,|\,A_b\,|\,A_c]$ by taking two columns of $A_c$ with $I_2$ on rows $r,s$ and 1's on row $p$ and then  one column of $A_b$ with 1 on row $r$ and so 0's on rows $s,p$ and one column of $A_a$ with a 1 on row $s$ and so 0's on rows $r,p$.
 
   We wish to assert that  $\ncols{V(1,{\even})}=\sum_{i\in W(1,{\footnotesize \even}) }|Y_i|\le 2m$. We consider the set system ${\cal Y}$ with sets $Y_i$ for $i\in W(1,\even)$. We set $I=W(1,\even)$ and appeal to \lrf{inductionforY} below to obtain $\sum_{i\in W(1,{\footnotesize \even})}|Y_i|\le 2m$.

Thus we have shown $|V(1,\even)|\le 2m$. The same will hold for $V(1,\odd)$ since we never use the parity in our argument other than to ensure for $a,b\in W(j,\odd)$ that $|a-b|\ge 2$. Also the same holds for $V(2,\even)$, $V(2,\odd)$ by taking (0,1)-complements. Thus $|V(1,\odd)|\le 2m$, $|V(2,\even)|\le 2m$ and  $|V(2,\odd)|\le 2m$.  Now this has included all columns of $A$ with the exception of $A_k$ for which $|A_k|\le 2$ and hence for at most 2m columns. We now conclude that $A$ has at most $10m$ columns. \qed

\begin{lemma}\label{inductionforY} Let $I$ be an ordered set. Let ${\cal Y}=\{Y_i\,:\,i\in I\}$ be a system of distinct sets $Y_i\subseteq [m]$ for $i\in I$. Assume   $|Y_i\cap Y_j|\le 1$ for $i,j\in I$. Assume for all  triples $a,b,c\in I$ with $a<b<c$ with the property that $Y_c\cap Y_b=r$ and $Y_ c\cap Y_a=s$, must have $r=s$.  Then $\sum_{i\in I}|Y_i|\le 2m$.\end{lemma}

\proof We use induction on $m$ with the result being easy for $m=1$. Let $d$ be the maximum index in $I$.  
 
 Our first case is that $Y_d\cap Y_i=\emptyset$ for $i\in I\backslash d$. We form a new set family ${\cal Y}'={\cal Y}\backslash Y_d$, whose sets are indexed by $I'=I\backslash d$, and whose sets are contained in $[m]\backslash Y_d$. Thus $\sum_{i\in I\backslash d}|Y_i|\le 2(m-|Y_d|)$  (the case $Y_d=\emptyset$ also works this way) and so $\sum_{i\in I}|Y_i|\le 2(m-|Y_d|)+|Y_d|\le 2m$.  
 
 Our second case assumes $Y_d\cap Y_j=\{q\}$ for some $j\in I\backslash d$. 
Our properties yield $Y_d\cap Y_i=\emptyset$ or $Y_d\cap Y_i=\{q\}$ for all $i\in I\backslash d$.  Then form a new set family ${\cal Y}'={\cal Y}\backslash Y_d$, whose sets are indexed by $I'=I\backslash d$, and whose sets are contained in $[m]\backslash (Y_d\backslash q)$.  We verify that ${\cal Y}'$ has the desired properties on $m-|Y_d|+1$. We use  that $(Y_d\backslash q)\cap Y_i=\emptyset$ for $i\in I\backslash d$.  By induction  $\sum_{i\in I\backslash d}|Y_i|\le 2(m-|Y_d|+1)$ and so $\sum_{i\in I}|Y_i|\le 2(m-|Y_d|+1)+|Y_d|\le 2m$.
\qed

\vskip 10pt
The following result is needed to complete our knowledge of $\forb(m,{\cal F})$ for ${\cal F}\subset\{Q_1,Q_2,\ldots ,Q_9\}$. 

\begin{thm}We have that $\forb(m,\{Q_6,Q_7,Q_9\})$ is $O(m)$.\label{lasttriple}\end{thm}
\proof Let $A\in\Av(m,\{Q_6,Q_7,Q_9\})$. We proceed as above letting  $A_k$ be the columns of sum $k$ and apply \lrf{types}. We deduce that if $A_k$ is of type 1 then $\ncols{A_k}\le 2$ else $Q_6\prec A_k$. Similarly if $A_k$ is of type 2 then $\ncols{A_k}\le 2$ else $Q_7\prec A_k$. Thus $\ncols{A}\le 2m-2$. \qed

\begin{thm} Let ${\cal F}\subset\{Q_1,Q_2,\ldots ,Q_9\}$ with ${\cal F}\ne\emptyset$. If ${\cal F}\subseteq \{Q_1,Q_4,Q_6\}$ or if ${\cal F}\subseteq \{Q_2,Q_5,Q_7\}$ or if ${\cal F}\subseteq \{Q_6,Q_7,Q_8\}$
or if ${\cal F}\subseteq\{Q_6,Q_9\}$ or if ${\cal F}\subseteq\{Q_7,Q_9\}$ or if ${\cal F}=Q_3$ then $\forb(m,{\cal F})$ is $\Theta(m^2)$. In all other cases, $\forb(m,{\cal F})$ is $O(m)$. In those cases  $\forb(m,{\cal F})$ is $\Theta(m)$ or $\Theta(1)$ and \trf{classify} will determine the asymptotic growth rate of $\forb(m,{\cal F})$ as either $\Theta(m)$ or $\Theta(1)$ in those cases where $\forb(m,{\cal F})$  is $O(m)$.\label{anysubset}\end{thm}

\proof Given that $\forb(m,Q_i)$ is $\Theta(m^2)$ for $i\in [9]$, we need only demonstrate that $\forb(m,{\cal F})$ is $O(m)$ in the other cases.  We can use the results listed in \tbrf{tab:bounds} to identify all pairs $Q_i,Q_j$ with $\forb(m,\{Q_i,Q_j\})$ being $O(m)$. Consider this as yielding a graph on a vertex set $[9]$. Any subset $S\subset [9]$ which contains one of these pairs has $\forb(m,\cup_{i\in S}Q_i)$ being $O(m)$ by \remrf{FsubsetG}. For example, any superset of $\{Q_1,Q_4,Q_6\}$ contains a pair $Q_i,Q_j$ with $\forb(m,\{Q_i,Q_j\})$ being $O(m)$. In particular $\forb(m,\{Q_1,Q_2\})$, $\forb(m,\{Q_1,Q_3\})$, $\forb(m,\{Q_1,Q_5\})$, $\forb(m,\{Q_1,Q_7\})$, \break $\forb(m,\{Q_1,Q_8\})$, and $\forb(m,\{Q_1,Q_9\})$ are all $O(m)$. For example, any superset of  $\{Q_6,Q_9\}$ contains either contains a pair $Q_i,Q_j$ with $\forb(m,\{Q_i,Q_j\})$ being $O(m)$ or is a triple $\forb(m,\{Q_i,Q_j,Q_k\})$ with $\forb(m,\{Q_i,Q_j,Q_k\})$ being $O(m)$.   We have $\forb(m,\{Q_1,Q_9\})$, $\forb(m,\{Q_2,Q_6\})$, 
$\forb(m,\{Q_3,Q_6\})$, $\forb(m,\{Q_4,Q_9\})$, \break$\forb(m,\{Q_5,Q_6\})$, and $\forb(m,\{Q_8,Q_9\})$ are all $O(m)$. We have two exceptional pairs $\{Q_6,Q_7\}$ and $\{Q_7,Q_9\}$ but we have the  triple $\{Q_6,Q_7,Q_9\}$ for which \hfil\break$\forb(m,\{Q_6,Q_7,Q_9\})$ is $O(m)$ by \trf{lasttriple}.    
\qed

\vskip 10pt
We may summarize our investigations by saying the \corf{grand} when applied to a forbidden family predicts the correct asymptotic growth for a number of elementary cases. Perhaps the cases where \corf{grand} doesn't correctly predict the   asymptotic growth, such as  \trf{unusual}, are rare. It is premature to conjecture an analog of \corf{grand} for forbidden families.

\end{document}